\documentclass[10pt,leqno]{amsart} 

   \usepackage[centertags]{amsmath}
   \usepackage{amsfonts}
   \usepackage{amsmath}
   \usepackage{amssymb}
   \usepackage{amsthm}
   \usepackage{newlfont}
   \usepackage[latin1]{inputenc}
\usepackage{xcolor}

\allowdisplaybreaks[3]

\theoremstyle{plain}
\newtheorem{Theo}{Theorem}

\newtheorem{Lem}[Theo]{Lemma}
\newtheorem{TheoA}{Theorem}

\theoremstyle{definition}
\newtheorem{Def}{Definition}

\theoremstyle{remark}

\newcommand{\Ee}{\mathfrak{e}}

\newcommand{\lp}{\mathcal{L}\text{-}\mathcal{P}}

\newcommand{\RR}{\mathbb{R}}

\newcommand{\NN}{\mathbb{N}}

\numberwithin{equation}{section}

\newcommand\Hh{{\mathcal H}}
\newcommand\G{{\mathcal G}}

\newcommand\p{\mbox{$\mathfrak{p}$}}

\newcommand\Be{\mathfrak{b}}

\newcommand\sign{\operatorname{sign}}

\catcode`,\active

\catcode`\,12

\title[Asymptotic for zeros of Bell and Eulerian polynomials]{Asymptotic for the rightmost zeros of Bell and Eulerian polynomials}
\author{Antonio J. Dur\'an}
\address{Departamento de An\'a\-li\-sis Mate\-m\'a\-ti\-co and IMUS,
        Universidad de Sevilla,
        41080 Sevilla, Spain}
\email{duran@us.es}
   \date{}

   \thanks{This research was partially supported by PID2021-124332NB-C21
(Minis\-te\-rio de Cien\-cia e Inno\-va\-ci\'on and Feder Funds (European Union)), and
FQM-262 (Jun\-ta de Anda\-lu\-c\'ia).}



\begin{document}

\begin{abstract}
Write $\zeta_m(n)$, $1\le m\le n-1$, for the negative zeros of the $n$-th Bell polynomial, arranged in decreasing order.
In this paper, we prove the following asymptotic: for every positive integer $m$ we have
$$
\lim_{n\to \infty}\frac{\zeta_m(n)}{-m\left(\displaystyle\frac{m}{m+1}\right)^{n-1}}=1.
$$
The approach used to find this asymptotic applies to many other significant families of polynomials. In particular, analogous asymptotics are also proved for the negative rightmost zeros of Eulerian polynomials, $r$-Bell polynomials, linear combinations of $K$ consecutive Bell polynomials and many others.
\end{abstract}

 \maketitle

\section{Introduction}

In this paper we develop a method to produce  asymptotics for the rightmost zeros of some relevant families of polynomials which include the Bell and Eulerian polynomials.

The Stirling numbers of the second kind $S(n,j)$, $0\le j\le n$, are defined by
$$
S(n,j)=\sum_{i=0}^j\frac{(-1)^{j-i}i^n}{(j-i)!\, i!},\quad 0\le j\le n
$$
(so that
\begin{equation}\label{comp0}
x^n=\sum_{j=0}^nS(n,j)x(x-1)\cdots (x-j+1).)
\end{equation}
For the Stirling numbers we prefer the notation $S(n,j)$ used in \cite{bre}
instead of the also usual Karamata-Knuth notation $\left\{\begin{matrix} n\\j\end{matrix}\right\}$ \cite{kar,knu}.

The Bell polynomials $(\Be_n)_n$ are defined by
\begin{equation}\label{ben}
\Be_n(x)=\sum_{j=0}^nS(n,j)x^j,\quad n\ge 0.
\end{equation}
They can also be defined using the generating function
$$
e^{x(e^z-1)}=\sum_{n=0}^\infty \frac{\Be_n(x)}{n!}z^n,
$$
or the recurrence
\begin{equation}\label{rbel}
\Be_{n+1}(x)=x\left(1+\frac{d}{dx}\right)\Be_n(x),\quad n\ge 0, \quad \Be_0=1.
\end{equation}
And there are some other equivalent definitions (see \cite{boy} for a review of Bell polynomials).

Bell polynomials are sometimes called as Touchard or exponential polynomials and were studied by Ramanujan (see \cite[Chapter 3]{ber}) in his notebooks before they were introduced by Touchard \cite{tou} and Bell \cite{bel} (in fact, Bell and Touchard called them exponential polynomials because of its generating function). They are also call single-variable Bell polynomials to distinguish them from the $k$-variable Bell polynomials $\Be_n(x_1,\dots, x_k)$ (the single-variable case being a particular instance of the $k$-variable one:  $\Be_n(x)=\Be_n(x,\dots, x)$).

It was proved by Harper \cite{har} that the Bell polynomials have simple and non positive zeros, and that the zeros of $\Be_{n+1}$ strictly interlaces the zeros of $\Be_n$ (see Definition \ref{inter} below).
Bell polynomials have been mainly considered in combinatorics and number theory although they have applications in probability and some other areas.

We denote $\zeta_m(n)$, $1\le m\le n-1$, for the negative zeros of the $n$-th Bell polynomial, arranged in decreasing order (let us note that $\Be_n(0)=0$, $n\ge 1$).

Strong and weak asymptotic for the zeros of the Bell polynomials were found by Elbert (\cite{elb1} and \cite{elb2}). The following asymptotic of the leftmost zero $\zeta_{n-1}(n)$ was conjectured in \cite{MeCo}:
$$
\lim_{n\to \infty}\frac{\zeta_{n-1}(n)}{-n}=e.
$$
In \cite[p. 3]{kmp} it is explained how the conjecture follows from the Elbert results in \cite{elb1} and \cite{elb2}.

Our method produces the following asymptotic for the rightmost zeros of the Bell polynomials.

\begin{Theo}\label{onl}
For every positive integer $m$, we have
\begin{equation}\label{asi}
\lim_{n\to \infty}\frac{\zeta_m(n)}{-m\left(\displaystyle\frac{m}{m+1}\right)^{n-1}}=1.
\end{equation}
\end{Theo}

The following table illustrates the asymptotic with some numerical values. For $n=100$ and $1\le m\le 5$, we have

\begin{align*} &\hspace{1.3cm}\zeta_m(n) & &\hspace{5pt}-m\left(m/(m+1)\right)^{n-1} \\
&-1.577721810\times 10^{-30} & & -1.577721810\times 10^{-30}
\\ &-7.379005723\times 10^{-18}& &-7.378963279\times 10^{-18}
\\ &-1.284493401\times 10^{-12} & & -1.282880874\times 10^{-12}
\\ &-1.031939874\times10^ {-9}& & -1.018517988\times 10^{-9}
\\ &-7.547543309\times 10^{-8}& & -7.244804083\times 10^{-8}
\end{align*}
We have computed the first five rightmost zeros of $\Be_{100}(x)$ using Mathematica (with WorkingPrecision=200) and Maple (with Digits=200), obtaining the same result. So we are confident that all the digits included in the table above are correct (anywhere see my paper \cite{DPV}). This comment is also valid for the other tables below in this paper.

\medskip

We also consider the Eulerian polynomials. Eulerian numbers are defined by
$$
\left\langle \begin{matrix} n \\ j \end{matrix} \right\rangle =\sum_{i=0}^j(-1)^{i}\binom{n+1}{i}(j+1-i)^n,\quad 0\le j\le n
$$
(so that
\begin{equation}\label{comp1}
x^n=\sum_{j=0}^n\left\langle \begin{matrix} n \\ j \end{matrix} \right\rangle \binom{x+j}{n}.)
\end{equation}
Eulerian numbers where introduced by L. Euler \cite{eul} and have interest in combinatorics and some other areas.

For the Eulerian numbers we use the angle brackets notation which it is also used
in \cite{bre}.

The Eulerian polynomials $(\Ee_n)_n$ are defined by $\Ee_0=1$ and
\begin{equation}\label{eup}
\Ee_n(x)=\sum_{j=0}^{n-1}\left\langle \begin{matrix} n \\ j \end{matrix} \right\rangle x^j,\quad n\ge 1.
\end{equation}
At $x=-1$, they provide the value of the Riemann zeta function at the negative integers
$$
\Ee_n(-1)=(2^{n+1}-4^{n+1})\zeta(-n),\quad n\ge 0.
$$

The Eulerian polynomial $\Ee_n$ has $n-1$ simple and negative zeros (let us note that $\Ee_n$ has degree $n-1$), and the zeros of $\Ee_{n+1}$ strictly interlace the zeros of $\Ee_n$ \cite{law}.
If we write $\varsigma_m(n)$, $1\le m\le n-1$, for the negative zeros of the $n$-th Eulerian polynomial, arranged in decreasing order, we have the following result.

\begin{Theo}\label{onle}
For every positive integer $m$, we have
\begin{equation}\label{asie}
\lim_{n\to \infty}\frac{\varsigma_m(n)}{-\left(\displaystyle\frac{m}{m+1}\right)^{n}}=1.
\end{equation}
\end{Theo}

Since
$$
\left\langle \begin{matrix} n \\ j \end{matrix} \right\rangle=\left\langle \begin{matrix} n \\ n-j-1 \end{matrix} \right\rangle,
$$
the zeros of the Eulerian polynomials are invariant under the transformation $\varsigma\mapsto 1/\varsigma$. Hence Theorem \ref{onle} also provides the asymptotic for the leftmost zeros of the Eulerian polynomials.

The following table illustrates the asymptotic (\ref{asie}) with some numerical values. For $n=100$ and $1\le m\le 5$, we have

\begin{align*} &\hspace{1.3cm}\varsigma_m(n) & &\hspace{5pt}-\left(m/(m+1)\right)^{n-1} \\
&-7.888609052\times10^{-31} & & -7.888609052\times10^{-31}
\\ &-2.459673290\times10^{-18}& &-2.459654426\times10^{-18}
\\ &-3.212242943\times10^{-13} & & -3.207202185\times10^{-13}
\\ &-2.069305321\times10^{-10}& & -2.037035976\times10^{-10}
\\ &-1.266880108\times10^{-8}& & -1.207467347\times10^{-8}
\end{align*}

\medskip

The content of the paper is as follows.
In Section \ref{eul2} we explain our method to produce the asymptotics (\ref{asi}) and (\ref{asie}). Actually we will consider a broader class of polynomials: those
generated by applying multipliers on the following variant of Eulerian polynomials:
\begin{equation}\label{eup2i}
e_n(x)=x^n\Ee_n(1+1/x)=\sum_{j=0}^nj!\, S(n,j)x^j.
\end{equation}

\begin{Def}\label{mul} A sequence $T=(\lambda_n)_n$ of real numbers is called a multiplier sequence if, whenever the polynomial $p(x)=\sum_{j=0}^na_jx^j$ has only real zeros, the polynomial
$$
T(p(x))=\sum_{j=0}^n\lambda_j a_jx^j
$$
also has only real zeros.
\end{Def}

Examples of multipliers are the sequences $(1/n!)_n$  or $(a^{-n^2})_n$, with $a$ a real number, $a\ge 1$.

Given a multiplier $T=(\lambda_n)_n$ such that $\lambda_n>0$, $n\ge 0$, we consider in Section \ref{eul2} the sequence of polynomials
$$
\p_n^T(x)=T(e_n(x))=\sum_{j=0}^n\lambda_jj!\, S(n,j)x^j,\quad n\ge 0.
$$
It turns out that $\p_n^T$ has $n-1$ negative zeros (and a zero at $x=0$) and the zeros of $\p_{n+1}^T$ and $\p_n^T$ interlace.
If we write $\zeta_m^T(n)$ for the $m$-th negative zero of $\p_{n}^T$ (arranged in decreasing order), then we have the following asymptotic.

\begin{Theo}\label{onlti}
For every positive integer $m$, we have
\begin{equation}\label{asiti}
\lim_{n\to \infty}\frac{\zeta_m^T(n)}{\displaystyle-\frac{\lambda_m}{\lambda_{m+1}}\left(\frac{m}{m+1}\right)^{n}}=1.
\end{equation}
\end{Theo}

Theorem \ref{onl} is the case $\lambda_j=1/j!$ and Theorem \ref{onle} follows from the case $\lambda_j=1$.

In Section \ref{zer2} we consider the asymptotic for the rightmost zeros of the $r$-Bell polynomials (see Theorem \ref{onl2}).

In Section \ref{zer3}, we consider polynomials $p_n$ which are linear combinations of $K$ consecutive Bell polynomials (with constant coefficients). We prove that $p_n$ has always at least $n-K$ real zeros and some of them can be positive.
For $n$ big enough, we calculate the exact number of positive zeros which turns out to be constant. In this case, our asymptotic describes both the behaviour of all the positive zeros as well as the behaviour of the rightmost negative zeros (see Theorem \ref{onl3}).

\section{Polynomials generated by applying multipliers on (modified) Eulerian polynomials}\label{eul2}

As explained in the Introduction, Theorems \ref{onl} and \ref{onle} are consequences of a more general Theorem on
families of polynomials generated by using multipliers (see Definition \ref{mul}). We apply the multipliers to the following variant of the Eulerian polynomials

\begin{equation}\label{eup2}
e_n(x)=x^n\Ee_n(1+1/x),
\end{equation}
where $\Ee_n$ are the $n$-th Eulerian polynomial (\ref{eup}).
Using
$$
j!\, S(n,j)=\sum_{k=0}^n\left\langle \begin{matrix} n \\ k \end{matrix} \right\rangle \binom{k}{n-j}
$$
(see \cite[(6.39) p. 269]{GKP}), we have
\begin{equation}\label{eup3}
e_n(x)=\sum_{j=0}^nj!\, S(n,j)x^j.
\end{equation}
A simple computation using \cite[(7a) p. 10]{hir} shows that the polynomials $(e_n)_n$ satisfy
the recurrence
\begin{equation}\label{reul}
e_{n+1}(x)=x\left(1+(1+x)\frac{d}{dx}\right)e_n(x),\quad n\ge 0, \quad e_0=1.
\end{equation}

Taking into account the properties of the Eulerian polynomials, we deduce that the polynomial $e_n$ has $n-1$ simple and negative zeros (actually they are in the interval $(-1,0)$) and a zero at $x=0$, and the zeros of $e_{n+1}$ strictly interlace the zeros of $e_n$. Moreover, if
we write $\hat \varsigma_m(n)$, $1\le m\le n-1$, for the negative zeros of the polynomial $e_n$, arranged in decreasing order, we have that
\begin{equation}\label{sigs}
\hat \varsigma_m(n)=\frac{1}{\varsigma_{n-m}(n)-1}=\frac{\varsigma_{m}(n)}{1-\varsigma_{m}(n)}
\end{equation}
(because $1/\varsigma_{n-m}(n)=\varsigma_{m}(n)$), where $\varsigma_m(n)$ are the negative zeros of the Eulerian polynomial $\Ee_n$ arranged in decreasing order.

In order to apply multipliers to the sequence of polynomials $(e_n)_n$, We need some previous definitions and results.

\begin{Def}\label{def1} An entire function $f$ is said to be in the Laguerre-Pólya class if it can be expressed in the form
\begin{equation}\label{pspr}
f(z)=cz^me^{-az^2+bz}\prod_{j=1}^\infty \left(1-d_j z\right)e^{d_j z},
\end{equation}
where $a\ge 0$, $m\in \NN$, $b, c,d_j\in\RR$, $j\ge 1$, and $\sum_{j=1}^\infty d_j^2<+\infty$.
The Laguerre-Pólya class will be denoted by $\lp$.

If $f\in \lp$ has only negative zeros, then we will use the notation $f\in \lp (-\infty,0)$.

We say that an entire function $f$ is of type I (or first type) in the Laguerre-Pólya class, in short $f\in \lp I$, if
$f(z)$ or $f(-z)$ has a product representation of the form
$$
cz^me^{\alpha z}\prod_{j=1}^\infty \left(1+\zeta_jz\right),
$$
where $\alpha \ge 0$, $c\in \RR$, $m\in \NN$ and $\zeta_j> 0$, $\sum_{j=1}^\infty\zeta_j<\infty$.
\end{Def}

Multipliers (see Definition \ref{mul}) can be characterized using functions in the class $\lp I$ (see \cite{PS}, \cite{CrCs}).

\begin{TheoA}\label{PS2} A sequence $T=(\lambda_n)_n$ of real numbers is a multiplier if and only if the function
$$
f(z)=\sum_{n=0}^\infty \frac{\lambda_n}{n!}z^n
$$
is of type I in the Laguerre-Pólya class.
\end{TheoA}

We can generate a large class of multipliers using the following theorem proved by Laguerre (see \cite{CrCs}).

\begin{TheoA}\label{lag}
If $f\in \lp (-\infty,0)$, then the sequence $(f(n))_n$ is a multiplier.
\end{TheoA}

Take now a sequence of real polynomials $(p_n)_n$, $\deg p_n=n$, satisfying that
\begin{equation}\label{equ}
\begin{cases} \mbox{for $n\ge 0$, all the zeros of $p_n$ are real}, and &\\
\mbox{the zeros of $p_{n+1}$ interlace the zeros of $p_n$.}&\end{cases}
\end{equation}

We use in this paper the following definition of interlacing.

\begin{Def}\label{inter} We say that the zeros of the polynomials $p$ and $q$ with $\deg p=1+\deg q=n$ interlace if
$$
\alpha_1\le \beta_1\le \alpha_2\le \cdots \le \beta_{n-1}\le \alpha_n,
$$
where $\alpha_1\le  \cdots \le \alpha_n$ and $\beta_1\le \cdots \le \beta_{n-1}$ are the zeros of $p$ and $q$ respectively. If all the inequalities above are strict then we say that the zeros of $p$ and $q$ strictly interlace.
\end{Def}

We can generate new sequences of polynomials with the same properties using multipliers. Indeed, given a multiplier $T=(\lambda_n)_n$ with $\lambda_n\not=0$, $n\ge0$, then the sequence of polynomials $p_n^T(x)=T(p_n(x))$, $n\ge 0$, also satisfy (\ref{equ}).

The interlacing property is a direct consequence of the following theorem  (see \cite[Theorem 9]{bra}).

\begin{TheoA}\label{bra} Let $T$ be a linear operator acting in the linear space of real polynomials with the property that if $p$ has only real zeros then $T(p)$ has also only real zeros. If the zeros of $p$ and $q$ interlace, then the zeros of $T(p)$ and $T(q)$ also interlace.
\end{TheoA}

Given a multiplier $T=(\lambda_n)_n$ such that $\lambda_n>0$, $n\ge 0$, we define the sequence of polynomials
$$
\p_n^T(x)=T(e_n(x))=\sum_{j=0}^n\lambda_jj!\, S(n,j)x^j,\quad n\ge 0.
$$
According to the previous discussion, we have that $\p_n^T$ has $n$ real zeros and the zeros of $\p_{n+1}^T$ and $\p_n^T$ interlace. Since $\lambda_j>0$, we deduce that $\p_n^T$ has $n-1$ negative zeros and a zero at $x=0$.
If we write $\zeta_m^T(n)$, $1\le m\le n-1$, for the $m$-th negative zero of $\p_{n}^T$ (arranged in decreasing order), then we have the following asymptotic.

\begin{Theo}\label{onlt}
For every positive integer $m$, we have
\begin{equation}\label{asit}
\lim_{n\to \infty}\frac{\zeta_m^T(n)}{\displaystyle-\frac{\lambda_m}{\lambda_{m+1}}\left(\frac{m}{m+1}\right)^{n}}=1.
\end{equation}
\end{Theo}

Theorems \ref{onl} and \ref{onle} are particular cases of Theorem \ref{onlt}. Indeed, on the one hand, since obviously the operator $T$ defined by the sequence $\lambda_n=1$, $n\ge 0$, is a multiplier, and $\p_n^T=e_n$, Theorem \ref{onlt} gives the following asymptotic for the rightmost zeros $\hat \varsigma_m(n)$ of $(e_n)_n$:
$$
\lim_{n\to \infty}\frac{\hat\varsigma_m(n)}{\displaystyle-\left(\frac{m}{m+1}\right)^{n}}=1.
$$
Then the identity (\ref{sigs}) between the zeros of $e_n$ and $\Ee_n$ establishes the asymptotic (\ref{asie}) in Theorem \ref{onle}.

On the other hand, since $f(z)=1/\Gamma(z+1)\in \lp (-\infty,0)$, we have that the sequence $T=(1/n!)_n$ is a multiplier. Then
we get $\Be_n(x)=\p_n^T(x)$, which shows that Theorem \ref{onl} is a particular case of Theorem \ref{onlt}.

\medskip

In order to prove Theorem \ref{onlt}, we need to introduce some notation. Given two intervals $I=(a_I,b_I)$, $J=(a_J,b_J)$, we write $I\prec J$ if $b_I<a_J$, i.e., the interval $J$ is to the right of the interval $I$ and their intersection is empty. Hence, $I\prec J$ implies $I\cap J=\emptyset$.

Given the multiplier $T=(\lambda_n)_n$, we define
\begin{equation}\label{un}
u_n=\frac{\lambda_n}{\lambda_{n+1}},\quad n\ge 0.
\end{equation}
Since multipliers are log-concave sequences (see, for instance, \cite[(1.5), p. 242]{CrCs0}), that is, $\lambda_n^2\ge \lambda_{n-1}\lambda_{n+1}$, we have that $(u_n)_n$ is an increasing sequence.

We write $I_{m,n}\equiv I_{m,n}(\epsilon, T)$ for the interval
\begin{equation}\label{int}
I_{m,n}=\left(-(1+\epsilon)\frac{\lambda_m}{\lambda_{m+1}}\left(\frac{m}{m+1}\right)^{n},-(1-\epsilon)\frac{\lambda_m}{\lambda_{m+1}}\left(\frac{m}{m+1}\right)^{n}\right).
\end{equation}
A simple computation gives that, for $0<\epsilon<1/(2m+1)$,
\begin{equation}\label{int2}
I_{m,n}\prec I_{m,n+1}.
\end{equation}
Similarly, using that $(u_n)_n$ is increasing, we have that, for $0<\epsilon<1/(2(m+1)^2)$,
\begin{equation}\label{int3}
I_{m+1,n}\prec I_{m,n}.
\end{equation}

For a positive integer $m$, define the function
\begin{equation}\label{rho}
\rho_m(x)=(x+1)\left(\frac{m}{m+1}\right)^x,\quad x\ge 0.
\end{equation}

We will need the following lemma.

\begin{Lem}\label{lrho} For any fixed positive integer $m$, we have:
\begin{enumerate}
\item For $0\le j\le m-1$,  the sequence $\rho_m(j)$ is increasing;
\item $\displaystyle \rho_m(m-1)=\rho_m(m)=m\left(\frac{m}{m+1}\right)^{m-1}$;
\item for $m\le j$, $\rho_m(j)$ is decreasing;
\item there exists $N_m$ such that for $j\ge N_m\ge m+1$, then
$$
0<\rho_m(j)<\rho_m(N_m)<1,\quad 0<\left(1+\frac{1}{2(N_m+1)^2}\right)\rho_m(N_m)<1.
$$
\end{enumerate}
\end{Lem}

\begin{proof}
A simple computation shows that the function $\rho_m(x)$ attains a unique maximum at $x_M=-1-1/\log[m/(m+1)]$, increases when $x\in (0,x_M)$, decreases when $x\in (x_M,+\infty)$
and goes to $0$ when $x\to +\infty$. Using that
$$
\frac{x-1}{x}<\log(x)<x-1,\quad 0<x, \, x\not =1,
$$
we get $m-1<x_M<m$. From where the proof follows easily.
\end{proof}

Properties (1), (2) and (3) in Lemma \ref{lrho} are saying that the sequence $\rho_m(j)$, $j\ge 1$, is unimodal with a plateau of two points in $j=m-1,m$ (see \cite[section 7.1]{com} for the definition of unimodal sequences).

\medskip

For a positive integer $N$ and sequences $\tau_j(n)$, $0\le j\le N$, $n\in \NN$, such that
\begin{equation}\label{mth}
\mbox {$\tau_0(n)=1$ and $\lim_n\tau_j(n)=1$, $1\le j\le N$,}
\end{equation}
consider the polynomial defined by
\begin{equation}\label{lim}
q_{n,N}(x)=\sum_{j=0}^N\lambda_{j+1}(j+1)^{n}\tau_j(n)x^j.
\end{equation}
Let us note that, for $n$ big enough, $q_{n,N}$ has degree $N$.

The polynomial $q_{n,N}$ obviously depends on the sequences $\tau_j$ and the multiplier $T$, but in order to simplify the notation, we do not explicitly indicate this dependency.

\begin{Lem}\label{zqn}
Given $0<\epsilon<1/(2(N+1)^2)$, there exist a positive integer $n_0$ and a positive real number $\delta>0$
such that for $n\ge n_0$:
\begin{enumerate}
\item The polynomial $q_{n,N}$ (\ref{lim}) has $N$ simple and negative zeros.
\item Write $\zeta_{m,N}\equiv\zeta_{m,N}(n)$, $m=1,\dots, N$, for the zeros of $q_{n,N}$ arranged in decreasing order. Then
$\{i: \zeta_{i,N}\in I_{m,n}\}=\{m\}$, and hence, $q_{n,N}$ changes its sign in each interval $I_{m,n}$ (\ref{int}).
\item For $x\not\in \cup_{m=1}^NI_{m,n}$, we have
$$
|q_{n,N}(x)|\ge\delta.
$$
\end{enumerate}
(Let us remark that  $n_0$ and $\delta$ depends on $N,\epsilon$, $T$ and the sequences $\tau_j$, $1\le j\le N$, but not on $n$).
\end{Lem}

\begin{proof}
We first prove that for $n$ big enough (depending on $N,\epsilon$, $T$ and the sequences $\tau_j$, $1\le j\le N$), the polynomial $q_{n,N}$ has at least one zero in the interval $I_{1,n}$ (\ref{int}).
It is enough to prove that $q_{n,N}(x)$ changes its sign in $I_{1,n}$. Writing $a_\pm=1\pm \epsilon$, we have ($u_1$ is defined by (\ref{un}))
\begin{align*}
q_n(-u_1 a_\pm /2^{n})&=\sum_{j=0}^N\lambda_{j+1}(j+1)^{n}\tau_j(n)(-u_1 a_\pm /2^{n})^j\\
&=\sum_{j=0}^N(-1)^j\lambda_{j+1}(u_1 a_\pm )^j\tau_j(n)(\rho_1(j))^{n}\\
&=\lambda_1[1-\tau_1(n)a_\pm]+\sum_{j=2}^N(-1)^j\lambda_{j+1}(u_1 a_\pm )^j\tau_j(n)(\rho_1(j))^{n}.
\end{align*}
On the one hand, $1-\tau_1(n)a_\pm=1-\tau_1(n)\mp\tau_1(n)\epsilon$. On the other hand, Lemma \ref{lrho} says that $\rho_1(j)<\rho_1(1)=1$, $j\ge 2$. Hence,
from the hypothesis (\ref{mth}), we deduce that for $n$ big enough
$$
\sign [ q_n(-u_1 a_\pm /2^{n})]=\sign(1-\tau_1(n)a_\pm)=\mp 1.
$$
We also have that
\begin{equation}\label{cpc1}
| q_n(-u_1 a_\pm /2^{n})|\ge \lambda_1 \epsilon/2.
\end{equation}

\medskip
We next prove that for $n$ big enough, the polynomial $q_{n,N}$ has at least one zero in the interval $I_{2,n}$ (\ref{int}).
Proceeding as before, we get
\begin{align*}
q_n(-u_2 a_\pm (2/3)^{n})&=\sum_{j=0}^N\lambda_{j+1}(j+1)^{n}\tau_j(n)(-u_2 a_\pm (2/3)^{n})^j\\
&=-a_\pm\frac{\lambda_2^2}{\lambda_3}\left(\tau_1(n)-\tau_2(n)a_\pm\right)\left(\frac43\right)^{n}\\
&\quad \quad +\lambda_1+\sum_{j=3}^N(-1)^j\lambda_{j+1}(u_2 a_\pm )^j\tau_j(n)(\rho_2(j))^{n}.
\end{align*}
On the one hand,
$$
a_\pm\left(\tau_1(n)-\tau_2(n)a_\pm\right)
=(1\pm\epsilon)(\tau_1(n)-\tau_2(n))\mp \epsilon\left(1\pm\epsilon\right)\tau_2(n).
$$
On the other hand, Lemma \ref{lrho} says that $\rho_2(j)<\rho_2(2)=4/3$, $j\ge 3$. Hence, from the hypothesis (\ref{mth}), we deduce, for $n$ big enough,
$$
\sign[q_n(- u_2 a_\pm(2/3)^{n})]=-\sign\left[a_\pm\left(\tau_1(n)-\tau_2(n)a_\pm\right)\right]=\pm 1.
$$
We also have that
\begin{equation}\label{cpc2}
|q_n(-u_2 a_\pm (2/3)^{n})|\ge \frac{\lambda_2^2\epsilon}{2^2\lambda_3}.
\end{equation}
\medskip

In general, for $1\le m\le N$, proceeding as before, we get
\begin{align*}
q_n\left(-u_m a_\pm \left(\frac{m}{m+1}\right)^{n}\right)&=\sum_{j=0}^N\lambda_{j+1}(j+1)^{n}\tau_j(n)\left(-u_m a_\pm \left(\frac{m}{m+1}\right)^{n}\right)^j\\
&=\frac{(-1)^{m-1}\lambda_m^ma_\pm^{m-1}}{\lambda_{m+1}^{m-1}}\left(\tau_{m-1}(n)-\tau_m(n)a_\pm\right)\left(\rho_m(m)\right)^{n}\\
&\quad+\sum_{j=0,j\not =m-1,m}^N(-1)^j\lambda_{j+1}(u_m a_\pm )^j\tau_j(n)(\rho_m(j))^{n}.
\end{align*}
On the one hand, Lemma \ref{lrho} says that $\rho_m(j)<\rho_m(m)$, $0\le j\le N$, $j\not =m-1,m$.
Hence, taking into account (\ref{mth}), the dominant term (with respect to $n$) in the sum above for $q_n\left(-u_m a_\pm \left(\frac{m}{m+1}\right)^{n}\right)$ is the one corresponding to $\left(\rho_m(m)\right)^{n}$.

On the other hand, \begin{align*}
(-1)^{m-1}a_\pm^{m-1}&\left(\tau_{m-1}(n)-\tau_m(n)a_\pm\right)\\
&=(-1)^{m-1}a_\pm^{m-1}\left[(\tau_{m-1}(n)-\tau_m(n))\mp\epsilon\tau_m(n)\right],
\end{align*}
and so
$$
\lim _n(-1)^{m-1}a_\pm^{m-1}\left(\tau_{m-1}(n)-\tau_m(n)a_\pm\right)=\mp(-1)^{m-1}\epsilon a_\pm^{m-1} .
$$

Hence for $n$ big enough
$$
\sign\left[q_n\left(-u_m a_\pm \left(\frac{m}{m+1}\right)^{n}\right)\right]=\mp (-1)^{m-1}.
$$
We also have that
\begin{align}\nonumber
\left|q_n\left(-u_m a_\pm \left(\frac{m}{m+1}\right)^{n}\right)\right|&\ge \frac{\epsilon (1-\epsilon)^{m-1}\lambda_m^m}{2\lambda_{m+1}^{m-1}}\left(\rho_m(m)\right)^{n}
\\\label{cpc3}&\ge \frac{\lambda_m^m\epsilon}{2^{m}\lambda_{m+1}^{m-1}}
\end{align}
(since $\rho_m(m)\ge 1$, $m\ge 1$, and $1-\epsilon> 1/2$).

This proves that the polynomial $q_{n,N}$ has at least one zero in the interval $I_{m,n}$ (\ref{int}).

Since $0<\epsilon<1/(2(N+1)^2)$, the intervals $I_{m,n}$, $1\le m\le N$, are ordered from right to left and are disjoint (it is a consequence of (\ref{int3})). Hence $\{i: \zeta_{i,N}\in I_{m,n}\}=\{m\}$.
This also proves that the polynomial $q_{n,N}$ has $N$ simple and negative zeros for $n$ big enough.

Since $q_{n,N}$ is a polynomial of degree $N$, the estimates (\ref{cpc1}), (\ref{cpc2}) and (\ref{cpc3}) show that
$$
\inf_{x\not \in \cup_{m=1}^NI_{m,n}} |q_{n,N}(x)|\ge \min \left\{\frac{\lambda_m^m\epsilon}{2^{m}\lambda_{m+1}^{m-1}}:1\le m\le N\right\}=\delta>0
$$
(and $\delta$ does not depend on $n$).

\end{proof}

\medskip

The number of real zeros of the polynomial $q_{n,N}$ is, in general, smaller than $N$ for small values of $n$. For instance, for $\lambda_j=1/j!$, $n=1$ and $\tau_j(n)=1$, we have
$$
q_{1,N}(x) =\sum_{j=0}^N\frac{1}{j!}x^j,
$$
and  $q_{1,N}$ has no real zeros for $N$ even and one real zero for $N$ odd (this is a consequence of the identity
$$
q_{1,N}(x)=e^x\Gamma(n+1,x)/n!,
$$
where $\Gamma(n+1,x)$ denotes the incomplete Gamma function).

\medskip

We are now ready to prove Theorem \ref{onlt}.

\begin{proof}[Proof of Theorem \ref{onlt}]
Given a positive integer $m$, take $N_m\ge m$ as in the part (4) of the Lemma \ref{lrho}, so that
$$
\rho_{m}(j)\le \rho_{m}(N_m)<1,\quad N_m\le j.
$$
It is enough to prove that given $0<\epsilon<1/(2(N_m+1)^2)$ then for $n$ big enough $\zeta_m^T(n)\in I_{m,n}$ (\ref{int}).

Let us remark that we also have $(1+\epsilon)\rho_{m}(N_m)<1$ (because of the part (4) of the Lemma \ref{lrho}).

Writing
$$
q_{n,N_m}(x)=\sum_{j=0}^{N_m} \lambda_{j+1}(j+1)!\, S(n,j+1)x^j,\quad r_{n,N_m}(x)=\sum_{j=N_m+1}^{n-1} \lambda_{j+1}(j+1)!\, S(n,j+1)x^j,
$$
we have $\p_n^T(x)=x(q_{n,N_m}(x)+r_{n,N_m}(x))$.

The asymptotic $S(n,j)\sim j^n/j!$ (see \cite[p. 121]{bon}) for the Stirling numbers of the second kind shows that the polynomial $q_{n,N_m}$ satisfies the hypothesis of Lemma \ref{zqn}.

We next prove that
$$
|r_{n,N_m}(x)|\le \delta/2, \quad \mbox{for $|x|<(1+\epsilon)u_m(m/(m+1))^{n}$},
$$
where $\delta>0$ is the positive number provided by the Lemma \ref{zqn} and $u_m$ is given by (\ref{un}).

Using the estimates \cite[Theorem 3]{weg}
\begin{equation}\label{ttem}
S(n,j+1)\le \frac{(j+1)^{n}}{(j+1)!}
\end{equation}
which holds for $n$ big enough and all $0\le j\le n-1$,
we get, for $|x|<(1+\epsilon)u_m(m/(m+1))^{n}$ and $j\ge N_m$,
\begin{align*}
|\lambda_{j+1}(j+1)!\, S(n,j+1)x^j|&\le \lambda_{j+1}(j+1)^{n}[u_m(1+\epsilon)]^{j}\left(\frac{m}{m+1}\right)^{jn}\\
&=\left(\frac{(j+1)m^j}{(m+1)^j}\right)^{n}\lambda_{j+1}\left(\frac{\lambda_m}{\lambda_{m+1}}\right)^j(1+\epsilon)^{j}.
\end{align*}
Using that
$$
\frac{1}{u_j}=\frac{\lambda_{j+1}}{\lambda_j}
$$
is decreasing (see (\ref{un})) and the inequality of arithmetic and geometric means, we have, for $m\le N_m \le j\le n-1$,
\begin{align*}
\lambda_{j+1}^{1/j}&=\lambda_{m}^{1/j}\left(\frac{\lambda_{j+1}}{\lambda_j}\cdots \frac{\lambda_{m+1}}{\lambda_m}\right)^{1/j}=\lambda_{m}^{1/j}\left[\left(\frac{\lambda_{j+1}}{\lambda_j}\cdots \frac{\lambda_{m+1}}{\lambda_m}\right)^{1/(j-m+1)}\right]^{(j-m+1)/j}
\\ &\le \lambda_{m}^{1/j}\left[\frac{1}{j-m+1}\sum_{i=m}^j\frac{\lambda_{i+1}}{\lambda_i}\right]^{(j-m+1)/j}\le \lambda_{m}^{1/j}\left[\frac{\lambda_{m+1}}{\lambda_m}\right]^{(j-m+1)/j},
\end{align*}
and so
$$
\lambda_{j+1}\left(\frac{\lambda_m}{\lambda_{m+1}}\right)^j\le \lambda_{m}\left(\frac{\lambda_m}{\lambda_{m+1}}\right)^{m-1}.
$$
Hence
\begin{align*}
|\lambda_{j+1}(j+1)!\, S(n,j+1)x^j|&\le \left(\frac{(j+1)m^j}{(m+1)^j}\right)^{n}\lambda_{j+1}\left(\frac{\lambda_m}{\lambda_{m+1}}\right)^j(1+\epsilon)^{j}\\
&\le \lambda_{m}\left(\frac{\lambda_m}{\lambda_{m+1}}\right)^{m-1}\left(\rho_m(j)\right)^{n}(1+\epsilon)^{j}\\&
\le \lambda_{m}\left(\frac{\lambda_m}{\lambda_{m+1}}\right)^{m-1}\left((1+\epsilon)\rho_m(N_m)\right)^{n}.
\end{align*}
Writing
$$
\sigma=(1+\epsilon)\rho_m(N_m)<1,
$$
we have
for $|x|<(1+\epsilon)u_m(m/(m+1))^{n}$ and $n$ big enough,
$$
\left|r_{n,N_m}(x)\right|\le C\lambda_{m}\left(\frac{\lambda_m}{\lambda_{m+1}}\right)^{m-1}(n-N_m)\sigma^{n}\le \delta/2.
$$
Hence, for $|x|<(1+\epsilon)u_m(m/(m+1))^{n}$,
$$
|\p_n^T(x)|\ge |q_{n,N_m}(x)|-|r_{n,N_m}(x)|\ge |q_{n,N_m}(x)|-\delta/2.
$$
Write
$$
X=[-(1+\epsilon)u_m(m/(m+1))^{n},0]\setminus \cup_{m=1}^{N_m}I_{m,n}.
$$
On the one hand, using Lemma \ref{zqn}, we deduce that for $x\in X$,
$$
|\p^T_n(x)|\ge |q_{n,N_m}(x)|-\delta/2\ge \delta-\delta/2\ge \delta/2>0.
$$
And so $\p^T_n$ does not vanish in $X$.

On the other hand, since
$$
|q_{n,N_m}(x)|\ge\delta>\delta/2\ge |r_{n,N_m}(x)|,\quad \mbox{for $x\in X$},
$$
we deduce that $\sign(\p^T_n(x))=\sign(q_{n,N_m}(x))$, $x\in X$.
Hence it follows from Lemma \ref{zqn} that
$$
\p^T_n(-(1-\epsilon)u_m(m/(m+1))^{n})\p^T_n(-(1+\epsilon)u_m(m/(m+1))^{n})<0,
$$
and so $\p^T_n$ has an odd number of zeros in each interval $I_{m,n}$.

We next prove by induction on $k$, $1\le k\le m$, that actually $\p^T_n$ has exactly one zero in the interval $I_{k,n}$ and that this zero is $\zeta_k^T(n)$.

It turns out that since $0<\epsilon<1/(2(N_m+1)^2)$, we have that $I_{k,n}\prec I_{k,n+1}$ and $I_{k+1,n}\prec I_{k,n}$ (see (\ref{int2}) and (\ref{int3}), respectively). In particular, $I_{1,n}$ is the rightmost interval among the intervals $I_{k,n}$.

Take then $k=1$, we have that $\p^T_n(x)$ does not vanish in $(-(1-\epsilon)u_1/2^{n},0]\subset X$ and hence $\zeta_1^T(n)\in I_{1,n}$ (which it is the rightmost interval). If $\zeta_2^T(n)\in I_{1,n}$, since the negative zeros of $\p^T_n$ interlace the negative zeros of $\p^T_{n-1}$ then $\zeta_1^T(n-1)\in I_{1,n}$, but this is a contradiction because $\zeta_1^T(n-1)\in I_{1,n-1}$ and $I_{1,n}\cap I_{1,n-1}=\emptyset$.

Assume finally that the only zero of $\p^T_n$ in $I_{k,n}$ is $\zeta_k^T(n)$. Then since
$\p^T_n(x)$ does not vanish in $X$, we have that $\zeta_{k+1}^T(n)\in I_{k+1,n}$ (since $I_{k+1,n}\prec I_{k,n}$). If $\zeta_{k+2}^T(n)\in I_{k+1,n}$, since the negative zeros of $\p^T_n$ interlace the negative zeros of $\p^T_{n-1}$ then $\zeta_{k+1}^T(n-1)\in I_{k+1,n}$, but this is again a contradiction because $\zeta_{k+1}^T(n-1)\in I_{k+1,n-1}$ and $I_{k+1,n}\cap I_{k+1,n-1}=\emptyset$.
\end{proof}

We illustrate the asymptotic (\ref{asit}) with a couple of examples.

Applying  the multiplier $T=(1/n!)_n$ $(s+1)$-times to the polynomials $e_n$ (\ref{eup3}), Theorem \ref{onlt}  produces the asymptotic
\begin{equation}\label{asitx}
\lim_{n\to \infty}\frac{\zeta_m^T(n)}{\displaystyle-(m+1)^{s+1}\left(\frac{m}{m+1}\right)^{n}}=1
\end{equation}
for the  rightmost zeros of the polynomials
$$
\p_n^T(x)=\sum_{j=0}^n\frac{S(n,j)}{j!^s}x^j.
$$

The next table is for $s=1$, $n=100$ and $1\le m\le 5$:
\begin{align*} &\hspace{1.3cm}\zeta_{m}^T(n) & &\hspace{5pt}-(m+1)^2\left(m/(m+1)\right)^{n} \\
&-3.155443621\times10^{-30} & & -3.155443621\times10^{-30}
\\ &-2.213698534\times10^{-17} & &-2.213688984\times10^{-17}
\\ &-5.136682477\times10^{-12} & & -5.131523497\times10^{-12}
\\ &-5.148460916\times10^{-9}& & -5.092589941\times10^{-9}
\\ &-4.501541829\times10^{-7}& & -4.346882450\times10^{-7}
\end{align*}

\medskip

Since for $a>1$, $f(z)=a^{-z^2}/\Gamma(z+1)\in \lp (-\infty,0)$ (the class $\lp (-\infty,0)$ is closed for the product of functions), we have that the sequence $T=(a^{-n^2}/n!)_n$ is a multiplier. Hence, Theorem \ref{onlt} produces the asymptotic
\begin{equation}\label{asity}
\lim_{n\to \infty}\frac{\zeta_m^T(n)}{\displaystyle-ma^{2m+1}\left(\frac{m}{m+1}\right)^{n-1}}=1
\end{equation}
for the  rightmost zeros of the polynomials
$$
\p_n^T(x)=\sum_{j=0}^n\frac{S(n,j)}{a^{j^2}}x^j.
$$

The next table is for $a=2$, $n=100$ and $1\le m\le 5$:
\begin{align*} &\hspace{1.3cm}\zeta_{m}^T(n) & &\hspace{5pt}-m2^{2m+1}\left(m/(m+1)\right)^{n-1} \\
&-1.262177448\times10^{-29} & & -1.262177448\times10^{-29}
\\ &-2.361271644\times10^{-16}& &-2.361268249\times10^{-16}
\\ &-1.642602557\times10^{-10} & & -1.642087518\times10^{-10}
\\ &-5.231621042\times10^{-7}& & -5.214812099\times10^{-7}
\\ &-0.0001497935741& & -0.0001483735876
\end{align*}

\section{$r$-Bell polynomials}\label{zer2}

The $r$-Bell polynomials are defined from the $r$-Stirling numbers of the second kind.
Although $r$-Stirling numbers are usually defined for $r\in \NN$ (because of their combinatorics definition \cite{car,car2,bro}),
we define them here for $r\in \RR$. The $r$-Stirling numbers of the second kind,
$S_r(n+r,j+r)$, $r\in \RR$ and $0\le j\le n$, are then defined by
$$
(x+r)^n=\sum_{j=0}^nS_r(n+r,j+r)x(x-1)\cdots (x-j+1).
$$
The $r$-Bell polynomials $\Be_{n,r}$, $n\ge 0$, are defined by
$$
\Be_{n,r}(x)=\sum_{j=0}^nS_r(n+r,j+r)x^j
$$
(see, for instance, \cite{mez,MeCo,MeRa}).

The Stirling numbers of the second kind and the Bell polynomials are the case $r=0$ (see (\ref{comp0})).

For $r>0$, the $r$-Bell polynomials have simple and negative zeros, and the zeros of $\Be_{n+1,r}$ strictly interlace  the zeros of $\Be_{n,r}$ (see \cite{wil}, \cite{mez}).
We denote $\zeta_{m,r}(n)$, $1\le m\le n$, for the negative zeros of the $n$-th $r$-Bell polynomial, arranged in decreasing order.

\begin{Theo}\label{onl2}
For $r>0$ and every nonnegative integer $m$, we have
\begin{equation}\label{asis}
\lim_{n\to \infty}\frac{\zeta_{m,r}(n)}{-m\left(\displaystyle\frac{m+r-1}{m+r}\right)^{n}}=1.
\end{equation}
\end{Theo}

Since for $n\ge 1$, $S_r(n,0)=r^n$, the Bell polynomial $\Be_n$, $n\ge 1$, always vanishes at $x=0$, and so it has $n-1$ negative zeros. But $\Be_{n,r}(0)\not =0$, for $r>0$, and so it has $n$ negative zeros.
Hence for $r=0$ and $m=1$, the denominator in left hand side of (\ref{asis}) gives the zero at $x=0$ of the Bell polynomial $\Be_n$. Therefore for $r=0$ we have to change $m\mapsto m+1$ and, in doing that, we get the asymptotic (\ref{asi}) for the leftmost zeros of the Bell polynomials.

By performing some cosmetic changes, the proof of Theorem \ref{onl2} is just the same as that of Theorem \ref{onlt} in the previous section.

Indeed, the intervals $I_{m,n}$ (\ref{int}) have to be changed to
$$
I_{m,n}^r=\left(-(1+\epsilon)m\left(\frac{m+r-1}{m+r}\right)^{n},-(1-\epsilon)m\left(\frac{m+r-1}{m+r}\right)^{n}\right).
$$

In turn, the function $\rho_m$ (\ref{rho}) has to be changed to
\begin{equation}\label{rhor}
\rho_{m,r}(x)=(x+r)\left(\frac{m+r-1}{m+r}\right)^x,\quad x\ge 0.
\end{equation}
The sequence $\rho_{m,r}(j)$, $j\in\NN$, is again unimodal with a plateau of two points in $j=m-1,m$, and
$$
\rho_{m,r}(m-1)=\rho_{m,r}(m)=\frac{(m+r-1)^m}{(m+r)^{m-1}}.
$$
So Lemma \ref{rho} has the corresponding $r$-version.

The polynomials $q_{n,N}$ have to be changed to
\begin{equation}\label{limr}
q_{n,N}^r(x)=\sum_{j=0}^N\frac{(j+r)^{n}}{j!}\tau_j(n)x^j.
\end{equation}

And finally, instead of the asymptotic $S(n,j)\sim j^n/j!$ and the estimation (\ref{ttem}) for the Stirling numbers of the second kind, we have to use
$$
1-\left(\frac{j+r-1}{j+r}\right)^nj<\frac{j!}{(j+r)^n}S_r(n+r,j+r)<1
$$
(see \cite[pp. 303]{mez}).

\medskip

We illustrate the asymptotic (\ref{asis}) with some numerical values.  The next table is for $r=1$, $n=100$ and $1\le m\le 5$:
\begin{align*} &\hspace{1.3cm}\zeta_{m,r}(n) & &\hspace{5pt}-m\left((m+r-1)/(m+r)\right)^{n} \\
&-7.888609052\times 10^{-31} & & -7.888609052\times 10^{-31}
\\ &-4.919334005\times 10^{-18}& &-4.919308853\times 10^{-18}
\\ &-9.632945556\times 10^{-13} & & -9.621606556\times 10^{-13}
\\ &-8.251338404\times 10^{-10}& & -8.148143905\times 10^{-10}
\\ &-6.282927695\times 10^{-8}& & -6.037336736\times 10^{-8}
\end{align*}

\medskip

The following one for $r=5/4$, $n=100$ and $1\le m\le 5$:
\begin{align*} &\hspace{1.3cm}\zeta_{m,r}(n) & &\hspace{5pt}-m\left((m+r-1)/(m+r)\right)^{n} \\
&-2.969952406\times 10^{-26} & & -2.969952405\times 10^{-26}
\\ &-2.142691182\times 10^{-16}& &-2.142622735\times 10^{-16}
\\ &-6.724335278\times 10^{-12} & & -6.707556623\times 10^{-12}
\\ &-2.708998681\times 10^{-9}& & -2.660863864\times 10^{-9}
\\ &-1.405094382\times 10^{-7}& & -1.339363965\times 10^{-7}
\end{align*}

\section{Linear combinations of $K$ consecutive Bell polynomials}\label{zer3}

In this section, we consider families of polynomials which are linear combinations of a fixed number $K$ of consecutive Bell polynomials with constant coefficients.
Hence, we take real numbers $\gamma_i$, $0\le i\le K$, with $\gamma_0=1$, and define
\begin{equation}\label{fcb}
p_n(x)=\sum_{j=0}^K\gamma_{j}\Be_{n-j}(x),
\end{equation}
where we take $\Be_{i}=0$ for $i<0$.

From (\ref{rbel}) it follows that, for $n\ge K$,
$$
p_{n+1}(x)=x\left(1+\frac{d}{dx}\right)p_n(x).
$$
This fact has important consequences on the zeros of the polynomials $p_n$, $n\ge K$. In order to show that, we need the following definition.

\begin{Def} We say that a linear operator $T$ acting in the linear space of real polynomials with $\deg(T(p))>\deg (p)$ is a real zero increasing operator if, for all polynomial $p$, the number of real zeros of $T(p)$ is greater than the number of real zeros of $p$.
\end{Def}

We then have the following.

\begin{Lem}\label{uyl}
The operator $D=x\left(1+\frac{d}{dx}\right)$ is a real zero increasing operator. Hence the number of real zeros of the polynomial $p_n$, $n\ge K$, is at least $n-K+1$ if $K$ is odd and $n-K$ if $K$ is even.
Moreover, if $\zeta<\xi$ are two simple, real and consecutive zeros of $p$ of equal sign, then the polynomial $D(p)$ has to vanish in $(\zeta,\xi)$, and if $\zeta$ is a zero of $p$ of multiplicity $l>1$, then it is a zero of $D(p)$ of multiplicity $l-1$.
\end{Lem}

\begin{proof}
Let $p$ be a polynomial. Assume first that $p$ has even degree, hence the number of real zeros has to be even (we only consider polynomials with real coefficients): say then that $p$ has exactly $2k$ real zeros. Write $r(x)=p(x)e^x$. It is clear that $r$ has also exactly $2k$ real zeros and hence, because of Rolle Theorem, $r'(x)=(p(x)+p'(x))e^x$ has at least $2k-1$ real zeros. So $p(x)+p'(x)$ has at least $2k-1$ real zeros. Since the degree of $p(x)+p'(x)$ is even, we deduce that it has at least $2k$ real zeros, and so $x(p(x)+p'(x))$ has at least $2k+1$ real zeros. If the degree of $p$ is odd, we can proceed similarly.
The rest of the lemma follows easily because $xr'(x)=e^xD(p)(x)$.

\end{proof}

The polynomials $p_n$ (\ref{fcb}) can have positive zeros. We first show that there exists $j_0$ such that for all $n$ the number of positive zeros of the polynomial $p_n$ is bounded by $j_0$.
Indeed, define
\begin{equation}\label{cag}
a_{n,j}=\sum_{i=0}^K\gamma_iS(n-i,j),\quad 0\le j\le n,
\end{equation}
so that
$$
p_n(x)=\sum_{j=0}^na_{n,j}x^j.
$$
Using the asymptotic $S(n,j)\sim j^n/j!$, we get
$$
a_{n,j}\sim \sum_{i=0}^K\gamma_i\frac{j^{n-i}}{j!}=\frac{j^{n}}{j!}\left(1+\sum_{i=1}^K\frac{\gamma_i}{j^i}\right).
$$
Since for $j\to \infty$ we have $\sum_{i=1}^K \gamma_i/j^i\to 0$, we deduce that there exist $j_0$ such that $a_{n,j}>0$, for $j\ge j_0$. Hence $p_n$ has at most $j_0$ positive zeros.

Lemma \ref{uyl} them says that the number of negative zeros $n_-$ of $p_n$ is at least $n-K-j_0+1$. Write $\xi_j^-(n)$, $1\le j\le n_-$, for the $n_-$ negative zeros of $p_n$ arranged in decreasing order (and taking into account their multiplicity).

Using our asymptotic analysis, we actually will calculate exactly the exact number of positive zeros of $p_n$ for $n$ big enough. To do that we have to consider
the polynomial
\begin{equation}\label{pop}
P(x)=\sum_{j=0}^K\gamma_{j}x^{K-j}.
\end{equation}
Assuming that
\begin{equation}\label{2a}
P(m)\not =0, \quad \mbox{for all positive integer $m$},
\end{equation}
we define the set
$$
\mathcal H=\{l: \mbox{$l$ is a positive integer and $P(l)P(l+1)<0$}\}.
$$
Denote by $s_\Hh$ the number of elements of $\Hh$, so that
\begin{equation}\label{3a}
\Hh=\{h_i:i=1,\dots ,s_\Hh \},
\end{equation}
where $h_i<h_j$ if $i<j$. Let us note that since $P$ has degree $K$, then $s_\Hh\le K$.
Define finally
\begin{equation}\label{4a}
\G=\{1,2,3,\dots\}\setminus \mathcal H=\{g_i:1\le i\},
\end{equation}
where $g_i<g_j$ if $i<j$.

The following Theorem describes the asymptotic behavior of the positive zeros of $p_n$, as well as that of the rightmost negative zeros.

\begin{Theo}\label{onl3} Let $\gamma_i$, $1\le i\le K$, be real numbers with $\gamma_0=1$.  Assume that the polynomial $P$ (\ref{pop}) satisfies (\ref{2a}). Then, there exists $n_0$ such that for $n\ge n_0$ the polynomial $p_n$ (\ref{fcb}) has exactly $s_\Hh$ positive zeros and they are simple. Write then $\xi_m^+(n)$, $1\le m\le s_\Hh$, for the  positive zeros of the polynomial $p_n$
arranged in increasing order.
Then for any given positive integer $m$, we have
\begin{align}\label{asis3}
\lim_{n\to \infty}\displaystyle\frac{\xi_{m}^+(n)}{-\frac{h_mP(h_m)}{P(h_m+1)}\left(\frac{h_m}{h_m+1}\right)^{n-K-1}}&=1,\quad 1\le m\le s_\Hh,\\\label{asis4}
\lim_{n\to \infty}\displaystyle\frac{\xi_{m}^-(n)}{-\frac{g_mP(g_m)}{P(g_m+1)}\left(\frac{g_m}{g_m+1}\right)^{n-K-1}}&=1,\quad 1\le m.
\end{align}
Moreover, for any positive integer $m$, $n_0$ can also be taken so that for $n\ge n_0$ the $m$ rightmost negative zeros of the polynomial $p_n$ are simple.
\end{Theo}

Again, the proof of Theorem \ref{onl3} follows by performing some changes in the proof of Theorem \ref{onlt}. However these changes need further explanation than for the case of $r$-Bell polynomials.

First of all the intervals $I_{m,n}$ (\ref{int}) have to be changed to
$$
I_{m,n}^K=\left(-\frac{(1+\epsilon)mP(m)}{P(m+1)}\left(\frac{m}{m+1}\right)^{n-K-1},-\frac{(1-\epsilon)mP(m)}{P(m+1)}\left(\frac{m}{m+1}\right)^{n-K-1}\right),
$$
for $m\in \G$ (\ref{4a}), and
$$
I_{m,n}^K=\left(-\frac{(1-\epsilon)mP(m)}{P(m+1)}\left(\frac{m}{m+1}\right)^{n-K-1},-\frac{(1+\epsilon)mP(m)}{P(m+1)}\left(\frac{m}{m+1}\right)^{n-K-1}\right),
$$
for $m\in \Hh$ (\ref{3a}). Let us remark that for $m\in\G$, $I_{m,n}^K\subset (-\infty,0)$ and for $m\in\Hh$, $I_{m,n}^K\subset (0,+\infty)$.

The monotonicity properties (\ref{int2}) and (\ref{int3}) for the intervals slightly change depending on whether the intervals are on the positive part of the real line or in the negative part. Indeed, for $0<\epsilon <1/(2m+1)$, a simple computation shows that
$$
I_{m,n}^K\prec I_{m,n+1}^K,\quad m\in \G,
$$
but
$$
I_{m,n+1}^K\prec I_{m,n}^K,\quad m\in \Hh.
$$
As for the monotonicity property (\ref{int3}), we have that for a positive integer $N$, there exists $n_N$ such that, for all $0<\epsilon<1/3$, then
\begin{equation}\label{ints1}
I_{m+1,n}^K\prec I_{m,n}^K,\quad \hbox{$m,m+1\in \G$, $n\ge n_N$ and $m\le N$,}
\end{equation}
but
\begin{equation}\label{ints2}
I_{m,n}^K\prec I_{m+1,n}^K,\quad \hbox{$m,m+1\in \Hh$, $n\ge n_N$ and $m\le N$.}
\end{equation}
Indeed, take $m,m+1\in\G$. Then, (\ref{ints1}) is equivalent to prove
$$
-(1-\epsilon)(m+1)u_{m+1}\left(\frac{m+1}{m+2}\right)^{n-K-1}<-(1+\epsilon)mu_m\left(\frac{m}{m+1}\right)^{n-K-1},
$$
where, to simplify the notation, we have written $u_m=P(m)/P(m+1)$. A simple computation gives that the previous inequality is equivalent to
\begin{equation}\label{ints3}
\epsilon<\frac{1-\alpha_{m,n}}{1+\alpha_{m,n}},
\end{equation}
where
$$
0<\alpha_{m,n}=\frac{mu_m}{(m+1)u_{m+1}}\left(\frac{m(m+2)}{(m+1)^2}\right)^{n-K-1}.
$$
Write $\eta=\sup\{u_m/u_{m+1}=P(m)P(m+2)/P^2(m+1):m\in \G\}$, so that
$$
0<\alpha_{m,n}<\eta \left(\frac{m(m+2)}{(m+1)^2}\right)^{n-K-1}.
$$
Since $(m(m+2)/(m+1)^2)^{n-K-1}<1$ and it increases with $m$, it is easy to deduce that there exists $n_N$ such that for $n\ge n_N$ and $m\le N$ then
$\alpha_{m,n}<1/2$. And so (\ref{ints3}) holds for all $0<\epsilon <1/3$.

The proof of (\ref{ints2}) is similar.

The function $\rho_m$ (\ref{rho}) does not need to be changed.

But the polynomials $q_{n,N}$ have to be changed to
\begin{equation}\label{lims}
q_{n,N}^K(x)=\sum_{j=0}^N P(j+1)\frac{(j+1)^{n-K-1}}{j!}\tau_j(n)x^j.
\end{equation}

With these changes, we have the following lemma. The proof is similar to that of Lemma \ref{zqn} and it is omitted.

\begin{Lem}
Let $N$ and $\epsilon$ be a positive integer and a real number, respectively, such that  $0<\epsilon<1/(2N+1)$.
Then there exists a positive integer $n_0$ and a positive real number $\delta>0$ such that for $n\ge n_0$:
\begin{enumerate}
\item The polynomial $q_{n,N}^K$ (\ref{lims}) has $s_\Hh^N$ simple and positive zeros and $N-s_\Hh^N$ simple and negative zeros, where $s_\Hh^N$ is the number of elements of the set $\{i\in \Hh: i\le N\}$.
\item For $N\ge \max \Hh$, write $\xi_{m,N}^+(n)$, $m=1,\dots, s_\Hh$, for the positive zeros of $q_{n,N}^K$ arranged in increasing order, and $\xi_{m,N}^-(n)$, $m=1,\cdots, N-s_\Hh$, for the negative zeros of $q_{n,N}^K$ arranged in decreasing order. Then for $h_m\in \Hh$, $\{i: \zeta_{i,N}^+\in I_{h_m,n}^K\}=\{m\}$, and for $g_m\in \G$, $\{i: \zeta_{i,N}^-\in I_{g_m,n}^K\}=\{m\}$. Hence, $q_{n,N}^K$ changes its sign in each interval $I_{m,n}^K$.
\item For $x\not\in \cup_{m=1}^NI_{m,n}^K$, we have
$$
|q_{n,N}^K(x)|\ge\delta.
$$
\end{enumerate}
(Let us note that  $n_0$ and $\delta$ depends on $N,\epsilon$ and the sequences $\tau_j$, $1\le j\le N$, but not on $n$).
\end{Lem}

Finally, instead of the asymptotic $S(n,j)\sim j^n/j!$ for the Stirling numbers of the second kind, we have to use
$$
a_{n,j}\sim P(j)\frac{j^{n-K}}{j!}
$$
(which easily follows from (\ref{cag})).

Then Theorem \ref{onl3} can be proved proceeding as in the proof of Theorem \ref{onlt}.

\medskip

The following table illustrates the asymptotic (\ref{asis3}) with some numerical values.
We take
$$
\gamma_1=1,\quad  \gamma_2=-\frac{17}2, \quad \gamma_3=\frac{87}4, \quad \gamma_4=-\frac{135}{8},
$$
so that
$$
P(x)=(x-3/2)(x-5/2)(x-9/2).
$$
This gives $\Hh=\{1,2,4\}$ and $\G=\{3,5,6,7,\dots\}$.
Using Maple, we have checked that the polynomial $p_n$, $n\ge 14$, has $s_\Hh=3$ positive zeros and $n-4$ negative zeros. Since $P$ satisfies the assumption (\ref{2a}), Theorem \ref{onl3} gives the asymptotic for the three positive zeros and the rightmost negative zeros of the polynomial $p_n$ (\ref{fcb}). For $n=100$ and $1\le m\le 3$, we have for the three positive zeros:

\begin{align*} &\hspace{1.3cm}\xi_m^+(n) &  &\hspace{0pt}-P(h_m)h_m^{n-2}/(P(h_m+1)(h_m+1)^{n-3}) \\
&5.301145283\times10^{-29}& &5.301145283\times10^{-29}
\\ &1.383545123\times10^{-17}& &1.383555614\times10^{-17}
\\ &8.250493853\times10^{-10} & & 8.525541195\times10^{-10}
\end{align*}

And for $n=100$ and $1\le m\le 3$, we have, for the three rightmost negative zeros,
\begin{align*} &\hspace{1.3cm}\xi_m^-(n) & &\hspace{0pt}-P(g_m)g_m^{n-2}/(P(g_m+1)(g_m+1)^{n-3}) \\
&-1.820667740\times10^{-12}& & -1.824541687\times10^{-12}
\\ &-2.469314309\times10^{-8}& &-2.318337306\times10^{-8}
\\ &-9.256545664\times10^{-7} & & -8.572668913\times10^{-7}
\end{align*}

\bigskip

We point out that in the forthcoming paper \cite{dur}, using the so-called generalized Bell polynomials, and assuming that the polynomial $P$ (\ref{pop}) has only real zeros and (\ref{2a}), we will prove that the polynomial $p_n$ (\ref{fcb}) has only real zeros for $n$ big enough.

\bigskip

Actually, Theorem \ref{onl3} can be extended to linear combination of $K$ consecutive (modified) Eulerian polynomials $(e_n)_n$ (\ref{eup3}).
Indeed for real numbers $\gamma_i$, $0\le i\le K$, with $\gamma_0=1$, we define
\begin{equation}\label{fce}
p_n(x)=\sum_{j=0}^K\gamma_{j}e_{n-j}(x),
\end{equation}
where we take $e_{i}=0$ for $i<0$.

From (\ref{reul}) it follows that, for $n\ge K$,
$$
p_{n+1}(x)=x\left(1+(1+x)\frac{d}{dx}\right)p_n(x).
$$
Since the differential operator $D(p(x))=x\left(1+(1+x)\frac{d}{dx}\right)$ is a real zero increasing operator, we can prove a version of Lemma \ref{uyl} for it. Write then $\xi_j^-(n)$, $1\le j\le n_-$, for the $n_-$ negative zeros of $p_n$ (\ref{fce}) arranged in decreasing order (and taking into account their multiplicity). Proceeding as before, we have the following asymptotic for the positive zeros of $p_n$, as well as that of the rightmost negative zeros.

\begin{Theo}\label{onl4} Let $\gamma_i$, $1\le i\le K$, be real numbers with $\gamma_0=1$.  Assume that the polynomial $P$ (\ref{pop}) satisfies (\ref{2a}). Then, there exists $n_0$ such that for $n\ge n_0$ the polynomial $p_n$ (\ref{fce}) has exactly $s_\Hh$ positive zeros and they are simple. Write then $\xi_m^+(n)$, $1\le m\le s_\Hh$, for the  positive zeros of the polynomial $p_n$ arranged in increasing order.
Then for a given positive integer $m$, we have
\begin{align*}
\lim_{n\to \infty}\displaystyle\frac{\xi_{m}^+(n)}{-\frac{P(h_m)}{P(h_m+1)}\left(\frac{h_m}{h_m+1}\right)^{n-K-1}}&=1,\quad 1\le m\le s_\Hh,\\
\lim_{n\to \infty}\displaystyle\frac{\xi_{m}^-(n)}{-\frac{P(g_m)}{P(g_m+1)}\left(\frac{g_m}{g_m+1}\right)^{n-K-1}}&=1,\quad 1\le m.
\end{align*}
Moreover, for any positive integer $m$, $n_0$ can also be taken so that for $n\ge n_0$ the $m$ rightmost negative zeros of the polynomial $p_n$ are simple.
\end{Theo}

\medskip

\paragraph{Acknowledgements.}
The author would like to thank the referee for his/her very careful reading of this paper. His/her comments and suggestions  have for sure improved the paper.



\end{document}